\documentclass {amsart}

\usepackage{amsmath}
\usepackage{amssymb}
\usepackage{amsthm}
\usepackage{graphicx}

\input xy
\xyoption{all}

\newcommand{\g}[1]{\mathfrak{#1}}
\newcommand{\mc} [1]{\mathcal{#1}}

\theoremstyle{plain}
\newtheorem{theorem}{Theorem} %[section]
\newtheorem{lemma}[theorem]{Lemma}

\newtheorem{cor}[theorem]{Corollary}

\theoremstyle{definition}

\newtheorem{example}[theorem]{Example}

\newtheorem{remark}[theorem]{Remark}

\begin{document}

\title {Multiple Points in ${\bf P}^2$ and Degenerations to Elliptic Curves}

\author {Ivan Petrakiev}
\address{Department of Mathematics, Harvard University, Cambridge, MA 02138}
\email{petrak@math.harvard.edu}
\begin {abstract} We consider the problem of bounding the dimension of the linear system of curves in ${\bf P}^2$ of degree $d$ with prescribed multiplicities $m_1,\dots,m_n$ at $n$ general points (\cite{Hir1},\cite{Hir2}). We propose a new method, based on the work of Ciliberto and Miranda (\cite {CM1}, \cite {CM2}), by specializing the general points to an elliptic curve in ${\bf P}^2$. 
\end{abstract}

\thanks{The author was partially supported by the NSF Graduate Research Fellowship.}
\maketitle

\tableofcontents

\section{Introduction}

Let $P_1,\dots,P_n$ be a set of $n$ general points in ${\bf P}^2$. For any $n$-tuple of positive integers ${\bf m} = (m_1,\dots,m_n)$ consider the ``fat-point'' scheme
$$\Gamma^{({\bf m})} = \bigcup P_i^{(m_i)}.$$
Determining the dimension of the linear system $|\mc{I}_{\Gamma^{({\bf m})}} (d)|$ of $d$-ics in ${\bf P}^2$ passing through each point $P_i$ with multiplicity $m_i$ is an open problem of algebraic geometry. In the present work we propose a new technique that allows to give an upper bound on this dimension in some cases.

To setup the notation, let ${\bf P}'$ be the blow-up of ${\bf P}^2$ at $P_1,\dots, P_n$. Then, $Pic ({\bf P}') = {\bf Z} H \oplus {\bf Z}E_1 \oplus \dots \oplus {\bf Z} E_n,$ where
$H$ is the pull-back of a hyperplane in ${\bf P}^2$ and $E_i$ is the exceptional divisor at $P_i$. Define the line bundle
$$\mc {L}_{\bf m} \cong \mc{O}_{{\bf P}'} (dH - \sum_{i=1}^n m_i E_i),$$
so that $|\mc{L}_{\bf m}| \cong |\mc{I}_{\Gamma^{({\bf m})}} (d)|.$ In future, we will omit the subscript $_{\bf m}$ and will simply write $\mc{L}$. By Riemann-Roch, the expected dimension $v$ of $|\mc{L}|$ is
$$v = \chi(\mc{L}) - 1 = \frac {d(d+3)} 2 - \sum_{i=1}^n \frac {m_i (m_i+1)}2.$$ 
We say that the linear system $|{\mc L}|$ is {\it special} if both cohomology groups $H^0({\mc L})$ and $H^1 ({\mc L})$ are nontrivial. We say that $|{\mc L}|$ is {\it homogeneous} if all multiplicities $m_i$ are equal to some fixed $m$. We have the following: \\

{\bf Conjecture (Harbourne-Hirschowitz \cite{CM1}, \cite{Hir2}).} {\it The linear system $|{\mc L}|$ is special if and only if $Bs(|{\mc L}|)$ contains a (-1)-curve $D$ with multiplicity at least two.} \\

In the homogeneous case, the conjecture would imply that there are no special linear systems with $n \geq 9$ (see \cite{CM2}). \\

Recently, Ciliberto and Miranda \cite{CM2} verified the Harbourne-Hirschowitz conjecture for all homogeneous linear systems $|\mc L|$ with $m \leq 12$. The basic idea is to specialize some of the general points to a line and study the degeneration of the linear system $|\mc L|$. 

In the present work, we propose to specialize the general points in ${\bf P}^2$ to a an elliptic curve instead of a line. 
%We begin by considering a motivating example in Section \ref{sec_motivation}. 
In Section \ref{sec_basic_construction}, we describe a degeneration of ${\bf P}^2$ into a union of two surfaces, namely a rational surface and an elliptic ruled surface. The basic construction, known as {\it the deformation to normal cone} (see \cite{F}), is similar to the one used by Ciliberto and Miranda in \cite{CM1}. 

In Section \ref{sec_main_r} we prove our main result (Theorem \ref{thm_main}), that gives a bound the dimension of $|\mc L|$ by the dimension of a (hopefully) simpler linear system in ${\bf P}^2$. 

Finally, in Section \ref{sec_limitations}, we give some applications of our result.

\begin{remark} The content of sections \ref{sec_basic_construction} and \ref{sec_main_r} generalizes to any smooth surface containing an elliptic curve, not just ${\bf P}^2$. We hope to find new interesting applications in future.
\end{remark}

\begin{remark} Specialization of multiple points to elliptic curves was also considered by Caporaso-Harris in unpublished notes \cite{CH}, where they used semi-stable reduction instead of deformation to normal cone. 
\end{remark}

\subsubsection*{Notation and Conventions} We work over an algebraically closed field of characteristic $0$. Recall some notation/terminology from \cite{Hart}. Let $C$ be a complete nonsingular curve. {\it A ruled surface} $S$ over $C$ is a nonsingular surface together with a ${\bf P}^1$-fibration $\pi : S \rightarrow C.$ {\it A minimal section} $C_0$ of $S$ is a section with minimal self-intersection. By a theorem of Atiyah (\cite{At}), $S$ is uniquely determined (upto a translation of $C$) by its {\it invariant} $e = -C_0^2$. For two divisors $Y$ and $Y'$ on $S$, $Y \sim Y'$ denotes rational equivalence and $Y \equiv Y'$ denotes numerical equivalence. Recall, that $Pic(S) = {\bf Z}C_0 \oplus Pic(C)$ and $Num(S) = {\bf Z} C_0 \oplus {\bf Z}f$, where $f$ is the class of a fiber. Thus, every divisor $Y$ on $S$ is rationally equivalent to some divisor $\mu C_0 + \g{b} f$, where $\g{b}$ is a divisor on $C$ and $\g{b} f := \pi^* (\g{b})$.

\subsubsection*{Acknowledgements}The author is grateful to J. Harris and S. Kleiman for enlightening discussions. The author thanks E. Cotterill for some useful remarks.

%\subsection {A Motivating Example}\label{sec_motivation}

%We consider a very concrete example 

%Consider the homogeneous linear system $|\mc{L}|$ of curves of degree $d=13$ passing through 10 general points $P_i \in {\bf P}^2$ with multiplicity $m=4$. The expected dimension of $|\mc{L}|$ is $v = 4$. We would like to show, that $|\mc{L}|$ is nonspecial as follows.

%Let $\mc{U} \in |\mc{L}|$ be any curve. Let's specialize the points $P_i$ to general points $P_i^* \in C$. Thus, we get a flat family of curves with a general member $\mc{U}$ and a special member $\mc{U}^*$. Since ${\mc U}^* \cdot E = 13 \times 3 - 10 \times 4 = -1$, Bezout's theorem implies, that the curve ${\mc U}^*$ contains $E$ as a component, with some multiplicity $\mu \geq 1$. In fact, since $({\mc U}^* - E) \cdot E = 0$ and the points $P_i^* \in C$ are general, we infer, that $\mu \geq 2$. Now, the ``residual'' curve $\mc{U}^*_0 = \mc{U}^* - 2E$ is of degree $13 - 2\times 3 = 7$ passing through the 10 points $P_i^*$ with multiplicity $4-2 = 2$. The dimension of the linear system of all such curves is 5, and so, by semicontinuity, $\dim |\mc{L}| \geq 5$. Unfortunately, this does not allow us to conclude, that $|\mc{L}|$ is nonspecial. 

\section {Basic construction}\label{sec_basic_construction}

%\begin{figure}
%\begin{center}
%\includegraphics[0in,0in][4in,6in]{fig1_1.eps}
%\caption{Blow-ups and a flip $(i \leq k)$.}
%\end{center}
%\end{figure}

Denote by $\Delta$ the affine line over the base field. The following lemma is motivated by the main construction in \cite{CM1} for degenerating ${\bf P}^2$. 

\begin{lemma}\label{lemma_construction} Fix positive integers $n \geq k \geq 10$. There exists a flat family of surfaces $X \rightarrow \Delta$ such that:

i) the general fiber $X_t$ is isomorphic to the blow-up of ${\bf P}^2$ at $n$ general points;

ii) the special fiber $X_0$ is the union of two components $S \cup {\bf P}'$ intersecting transversally along an elliptic curve $C$. Here, $S$ is an indecomposable ruled surface over $C$; the component ${\bf P}'$ is isomorphic to the blow-up of ${\bf P}^2$ at $n-k$ general points in ${\bf P}^2$ and $k$ general points on $C$.
\end{lemma}

\begin{proof} We describe the construction of $X$ in the following five steps 
%(see Fig. 1)
. The first two steps are just {\it the deformation to the normal cone} of an elliptic curve in ${\bf P}^2$ (\cite{F}). \\

{\it Step 1.} Let $X_1 = {\bf P}^2 \times \Delta$ be the trivial family of planes. Let ${\bf P}_0$ be the fiber of the projection map $X_1 \rightarrow \Delta$ at $t=0$. Fix a nonsingular elliptic curve $C \subset {\bf P}_0$. For any $i \leq n$, let $p_i: \Delta \rightarrow X_1$ be a section of the projection map $X_1 \rightarrow \Delta$. Denote by $P_i \in {\bf P}_0$ the image $p_i(0)$. We assume the following: $(i)$ $p_i$ is an embedding; $(ii)$ for $t$ general, $p_i(t)$ is a general point in ${\bf P}^2$; $(iii)$ for $i\leq k$, $P_i$ is a general point on $C$; $(iv)$ for $i > k$, $P_i$ is a general point in ${\bf P}_0$ and $(v)$ for $i \leq k$, the image of $p_i$ intersects ${\bf P}_0$ transversally at $P_i$. For example, we may take $p_i : \Delta \rightarrow X_1$ to be linear maps satisfying the above properties. \\

{\it Step 2.} Let $X_2 \rightarrow X_1$ be the blow-up of $X_1$ along $C$ and let $S_0$ be the exceptional divisor of the blow-up. Thus $S_0 = {\bf P} (\mc{O}_C \oplus \mc{O}_C (3))$ is an elliptic ruled surface with minimal section $C$. Hence, $S_0$ has invariant $e = 9$. Another way of seeing this is to use the triple-point formula (\cite{P}, Cor. 2.4.2), according to which $(C|_{S_0})^2 = -(C|_{{\bf P}_0})^2 = -9$. 

For any $i$, we lift the section $p_i$ to a section $\tilde p_i :\Delta \rightarrow X_2$. Denote by $\tilde P_i$ the image $\tilde p_i(0)$. Thus, $\tilde P_1,\dots, \tilde P_k$ are general points on $S_0$ and $P_{k+1}, \dots, P_n$ are general points on the proper transform of ${\bf P}_0$ (which we denote again by ${\bf P}_0$, abusing the notation).

For $i \leq k$, denote by $F_i$ the unique fiber of $S_0$ passing through the point $\tilde P_i$. From the exact sequence
$$0 \longrightarrow \mc{N}_{{F_i} / S_0} \longrightarrow \mc{N}_{{F_i} / X_2} \longrightarrow \mc{N}_{S_0 / X_2}|_{F_i} \longrightarrow 0$$
we find that $\mc{N}_{F_i / X_2} \cong \mc{O}_{F_i} \oplus \mc{O}_{F_i}(-1)$. \\

{\it Step 3.} We blow-up the images of $\tilde p_i$ simultaneously and denote by $X_3$ the resulting threefold. Let $E_i$ be the exceptional divisors of the blow-up. Now, $X_3$ is a family of surfaces over $\Delta$ such that: $(i)$ the general fiber $(X_3)_t$ is isomorphic to the blow-up of ${\bf P}^2$ at $n$ general points and $(ii)$ the special fiber is the union of two components $S_0^+ \cup {\bf P}^+_0$, where $S_0^+$ is the blow-up of $S_0$ at the $k$ points $\tilde P_1,\dots, \tilde P_k$ and ${\bf P}_0^+$ is the blow-up of ${\bf P}_0$ at the $n-k$ points $\tilde P_{k+1},\dots, \tilde P_n$. By a similar computation as in Step 2, we find that $\mc {N}_{F_i/X_3} \cong \mc{O}_{F_i} (-1) \oplus \mc{O}_{F_i} (-1)$. For any $t \in \Delta$, we denote by $E^{(t)}_i$ the restriction of $E_i$ to $(X_3)_t$. \\

{\it Step 4.} In the next two steps, we ``transfer'' the exceptional curves $F_i$ from $S_0^+$ to ${\bf P}_0^+$ by applying the {\it flip} transformation (see \cite{P}, Cor. 2.4.5 or \cite{Deb}, 6.20). 

We blow-up the curves $F_i$ simultaneously and denote by $X_4$ the resulting threefold. Denote by $W_i$ the exceptional divisors of the blow-up. Since $\mc {N}_{F_i/X_3} \cong \mc{O}_{F_i} (-1) \oplus \mc{O}_{F_i} (-1)$, we conclude that $W_i \cong {\bf P}^1 \times {\bf P}^1$. 

Let ${\bf P}_0^{++}$ be the proper image of ${\bf P}_0^+$ in $X_4$. Thus ${\bf P}_0^{++}$ is isomorphic to the blow-up of ${\bf P}^2$ at the $n$ points $P_1,\dots,P_n$. Let $D_1,\dots, D_k$ be the exceptional divisors of the blow-up ${\bf P}^{++}_0$ at $P_1,\dots, P_k$. Notice that $F_i$ and $D_i$ belong to different rulings of $W_i$. \\

{\it Step 5.} We contract all $W_i$ simultaneously along the ruling given by $F_i$ and denote by $X_5$ the resulting threefold. Let $S_0^{-}$ be the image of $S_0^+$ in $X_5$. Clearly, $S_0^{-}$ is obtained from $S_0$ (defined in Step 2) by applying $k$ {\it elementary transforms} at the points $\tilde P_1,\dots,\tilde P_k$ (\cite{Hart}, Example V.5.7.1). Since $S_0$ has invariant $e=9$, the points $\tilde P_1, \dots, \tilde P_k$ are general and $k \geq 10$, it follows that $S_0^-$ is an indecomposable ruled surface over $C$. Therefore, the threefold $X=X_5$ has the required properties (with $S = S^-_0$ and ${\bf P}' = {\bf P}^{++}_0$).
\end{proof}

\section {Main Result}\label{sec_main_r}

Fix positive integers $d,n,k$ and $m_1,\dots,m_n$ where $n \geq k \geq 10$. Let $X \rightarrow \Delta$ be the family of surfaces constructed in the previous section. For any $t$, we denote by $X_t$ the fiber of $X$ at $t$. Denote by $E_i$ the exceptional divisors on $X$ corresponding to the $n$ general points. For $t$ general, denote by $E^{(t)}_i$ the restriction $E_i|_{X_t}$. Thus, $E^{(t)}_i$ is an exceptional divisor of the blow-up $X_t \rightarrow {\bf P}^2$. For $t=0$ and $i \leq k$, the restriction $E_i |_{X_0}$ has two components $E^{(0)}_i$ and $D_i$ 
%(as in Fig. 1)
. Here, $E^{(0)}_i$ is just a fiber of the ruled surface $S$ and $D_i$ is an exceptional divisor of the blow-up ${\bf P}' \rightarrow {\bf P}^2$. If $i > k$, then $E_i |_{X_0}$ has only one component, denoted by $E_i^{(0)}$.

Consider the line bundle $\mc{L} = \mc{O}_X (dH - \sum_{i=1}^n m_i E_i)$ on the threefold $X$. In particular, $$\mc{L}|_{X_t} \cong \mc{O}_{X_t} (dH - \sum_{i=1}^n m_i E^{(t)}_i),$$ for $t$ general. At the special fiber, we have:
$$\mc{L}|_{{\bf P}'} \cong \mc{O}_{{\bf P}'} (dH - \sum_{i=1}^k m_i D_i - \sum_{i=k+1}^n m_i E^{(0)}_i)$$
and
$$\mc{L}|_S \cong \pi^* \mc{O}_C(dH - \sum_{i=1}^k m_i D_i) \cong \mc{O}_S (\g{b}f),$$
for a suitable divisor $\g{b}$ on $C$ (by construction, $\g{b}$ is general). Here $\pi$ denotes the projection $\pi : S \rightarrow C$.

For any integer $\mu$, consider the twist $\mc{L} (\mu) = \mc{L} \otimes \mc{O}_X (-\mu S)$. Since $\mc{O}_X(S + {\bf P}') \cong \mc{O}_{X} (X_t) \cong \mc{O}_X$, we conclude that $\mc{O}_X(-S) \cong \mc{O}_X ({\bf P}')$. Therefore,
$$\mc{L}(\mu)|_{{\bf P}'} \cong \mc{O}_{{\bf P}'} ((d-3\mu)H - \sum_{i=1}^k (m_i-\mu) D_i - \sum_{i=k+1}^n m_i E^{(0)}_i)
$$
and
$$\mc{L}(\mu)|_S \cong \mc{O}_S (\mu C + \g{b} f).$$

Notice, that $\mc{L}(\mu)|_{X_t} \cong \mc{L}|_{X_t}$ for $t \in \Delta$ general and any $\mu \in {\bf Z}$. Thus, we should think of $\mc{L}(\mu)|_{X_0}$ as a limit of the linear system $\mc{L}|_{X_t}$ as $t \rightarrow 0$. In particular, any choice $\mu$ leads to a possible limit (compare with the theory of limit linear series on curves, introduced by Eisenbud-Harris in \cite{EHLim}). \\

We are now in position to formulate the main result in this section.

\begin{theorem}\label{thm_main} Let $\mu$ be a positive integer such that $\chi(\mc{L}(\mu)|_{{\bf P}'}) \geq \chi(\mc{L}|_{X_t})$ for a general $t \in \Delta$. Then $h^0 (\mc{L}(\mu)|_{{\bf P}'}) \geq h^0 (\mc{L}|_{X_t})$.
\end{theorem}

The number $\mu$ should be interpreted as follows: let $\mc{U}$ be a curve in ${\bf P}^2$ passing through $n$ general points $P_1,\dots,P_n$ with multiplicity $m_1,\dots, m_n$. As we specialize the first $k$ of the points to an elliptic curve $C$ (in a general fashion), at least $\mu$ copies of $C$ must split-off from $\mc{U}$. \\

The following lemma plays an essential role in the proof of the theorem.

\begin{lemma} Let $S$ be an indecomposable ruled surface over an elliptic curve and let $C$ be a section of $S$. Let $D \sim \mu C + \g{b} f$ be an effective divisor on $S$, where $\mu > 0$ and $\g{b} \in Pic(C)$ is general. Then, $D$ is ample and $\chi(\mc{O}_S(D)) > 0$.
\end{lemma}

\begin{proof} Let $C_0$ be a minimal section of $S$ and let $e = -C_0^2$. We may write $D \sim \mu C_0 + \g{b}' f$, where $\g{b}' \in Pic(C)$ is general. The canonical divisor of $S$ is $K_S \equiv -2C_0 - ef$ and the arithmetic genus of $S$ is $p_a = -1$ (see \cite{Hart}, Ch. V.2). By Riemann-Roch, 
$$\chi(\mc{O}_S(D)) = \frac 1 2 D \cdot (D-K_S) + p_a + 1 = (\mu+1)(b' - \frac 1 2 {\mu e}),$$ where $b' = \deg{\g{b}}'$. Therefore, to show that $\chi(\mc{O}_S(D)) > 0,$ it sufficies to show that $$b' - \frac 1 2 {\mu e} > 0.$$

Since $S$ is indecomposable, $e = 0$ or $-1$ (\cite{Hart}, Thm. V.2.15). Suppose that $e = 0$. Since $C_0\cdot D = b'$ and $C_0$ is nef, we have $b' \geq 0$. In fact, $b' > 0$, because $\g{b}'$ is general (it suffices to assume that the line bundle $\mc{O}_{C_0}(D)$ is not a multiple twist of $\mc{O}_{C_0}(C_0)$). 

Suppose that $e = -1$. Then, it is well-known that $S$ contains a nonsingular elliptic curve $Y \equiv 2C_0 - f$ (see \cite{Deb}, p.24). Since $Y^2 = 0$, it follows that $Y$ is nef. Therefore, $Y \cdot D = 2b' + \mu \geq 0$. In fact, $2b' + \mu > 0$, because $\g{b}'$ is general (it suffices to assume that the line bundle $\mc{O}_Y(D)$ is not a multiple twist of $\mc{O}_Y(Y)$).

The fact that $D$ is ample follows from the description of the ample cone of $S$ (see \cite{Hart}, prop. V.2.20 and 2.21).

This completes the proof of the lemma.
\end{proof}

{\it Proof of the theorem.} It will be notationally more convenient to replace $\mu$ with $\mu+1$ in the statement of the theorem. In other words, given that $\chi(\mc{L}(\mu+1)|_{{\bf P}'}) \geq \chi(\mc{L}|_{X_t})$ we want to show that $h^0 (\mc{L}(\mu+1)|_{{\bf P}'}) \geq h^0(\mc{L}|_{X_t})$.

Consider the Mayer-Vietoris exact sequence on $X_0$:
$$0 \longrightarrow \mc{O}_{X_0} \longrightarrow \mc{O}_{{\bf P}'} \oplus \mc{O}_{S} \longrightarrow \mc{O}_C \longrightarrow 0.$$
We tensor the above sequence with $\mc{L}(\mu)$ and take cohomology:
$$
0 \longrightarrow H^0 (\mc{L}(\mu)|_{X_0}) \longrightarrow H^0(\mc{L}(\mu)|_{{\bf P}'}) \oplus H^0(\mc{L}(\mu)|_{S}) \mathop {\longrightarrow}^{f\oplus g} H^0 (\mc{L} (\mu)|_C)
$$
We have:
$$\chi(\mc{L}(\mu)|_{X_0}) + \chi(\mc{L}(\mu)|_C) = \chi (\mc{L}(\mu)|_{{\bf P}'}) + \chi (\mc{L}(\mu)|_S).$$
Consider the following exact sequence on ${\bf P}'$:
$$0 \longrightarrow \mc{O}_{{\bf P}'}(-C) \longrightarrow \mc{O}_{{\bf P}'} \longrightarrow\mc{O}_C \longrightarrow 0.$$
We tensor the above sequence with $\mc{L}(\mu)$ and take cohomology:
$$0 \longrightarrow H^0 (\mc{L}(\mu+1)|_{{\bf P}'}) \longrightarrow H^0 (\mc{L}(\mu)|_{{\bf P}'}) \mathop {\longrightarrow}^f H^0 (\mc{L} (\mu)|_C)$$
We have:
$$\chi(\mc{L}(\mu)|_{{\bf P}'}) = \chi (\mc{L}(\mu+1)|_{{\bf P}'}) + \chi (\mc{L}(\mu)|_C).$$
Adding the last two equalities gives:
$$\chi(\mc{L}(\mu)|_{X_0}) = \chi(\mc{L}(\mu+1)|_{{\bf P}'}) + \chi (\mc{L}(\mu)|_S). \eqno(*)$$
Since Euler characteristic is constant in flat families, we have $$\chi(\mc{L}(\mu)|_{X_0}) = \chi(\mc{L}|_{X_t} (\mu)) = \chi(\mc{L}|_{X_t}),$$ for $t$ general. Now, the assumption $\chi(\mc{L}(\mu+1)|_{{\bf P}'}) \geq \chi(\mc{L}|_{X_t}),$ together with (*), implies $\chi(\mc{L}(\mu)|_S) \leq 0.$ So, by the previous lemma, $$H^0(\mc{L}(\mu)|_S) = 0.$$ Now, from the last two exact sequences in cohomology, $$H^0 (\mc{L}(\mu)|_{X_0}) = \text{ker\ } f = H^0 (\mc{L}(\mu+1)|_{{\bf P}'}).$$
Finally, by semicontinuity,
$$h^0(\mc{L}|_{X_t}) = h^0(\mc{L}(\mu)|_{X_t}) \leq h^0(\mc{L}(\mu)|_{X_0}) = h^0 (\mc{L}(\mu+1)|_{{\bf P}'}).$$
This completes the proof. \hfill$\Box$\\

\section {Applications}\label{sec_limitations}

In this final section, we will use Theorem \ref{thm_main} to show that certain homogeneous linear systems in ${\bf P}^2$ are nonspecial. Also, we will give an example that exhibits a limitation of our theorem. 

Given data $(d,n,m)$, consider curves in ${\bf P}^2$ of degree $d$ passing through $n \geq 10$ general points with multiplicity $m$. For simplicity, we will specialize all $n$ points at once to an elliptic curve $C \subset {\bf P}^2$.

So, let $X \rightarrow \Delta$ and $\mc{L}$ be as before, with $k=n$. For any integer $\mu$, we have:
\begin{align*}
\chi(\mc{L}&(\mu)|_{{\bf P}'}) - \chi (\mc{L}|_{X_t}) \\
&= \frac {(d-3\mu)(d-3\mu+3)} 2 - n \frac {(m-\mu)(m-\mu+1)} 2 - \frac{d(d+3)} 2 + n \frac {m(m+1)} 2 \\
&= \frac 1 2 \mu (n-9-6d+2mn -\mu(n-9) ).
\end{align*}
In particular, $\chi(\mc{L}(\mu)|_{{\bf P}'}) \geq \chi (\mc{L}|_{X_t})$ if
$$0 \leq \mu \leq 1+\frac {2mn - 6d} {n-9}.$$
(Notice, that the right-hand side is just $1 + 2 (\mc{L}|_{X_t} \cdot K_{X_t}) / (-K_{X_t}^2)$ for $t \in \Delta$ general.) \\

Clearly, in order to get the most information from Theorem \ref{thm_main}, we should choose the greatest integral value of $\mu$, subject to the inequality above. The best scenario is achieved when the upper bound on $\mu$ is already an integer:

\begin{cor}\label{cor_cool} Let $(d,n,m)$ be as above, and assume, that $\mu = 1 + \frac {2mn - 6d} {n-9}$ is a positive integer. If $\mc{L}(\mu)|_{{\bf P}'}$ is nonspecial, then so is $\mc{L}|_{X_t}$, for $t$ general.
\end{cor}
\begin{proof} We have $\chi(\mc{L}(\mu)|_{{\bf P}'}) = \chi(\mc{L}|_{X_t})$ and $h^0 (\mc{L}(\mu)|_{{\bf P}'}) \geq h^0 (\mc {L} |_{X_t})$. Assuming that $h^0 (\mc{L}(\mu)|_{{\bf P}'}) > 0$, we have 
$$\chi (\mc{L}(\mu)|_{{\bf P}'}) = h^0 (\mc{L}(\mu)|_{{\bf P}'}) \geq h^0 (\mc{L}|_{X_t}) \geq \chi(\mc{L}|_{X_t}).$$
So, there is equiality everywhere. It follows, that $h^1 (\mc{L}|_{X_t}) = 0$.
\end{proof}

We proceed with some examples.

\begin{example} Consider the linear system corresponding to the data $(d,n,m) = (13,10,4)$, with expected dimension $v = \chi(\mc{L}|_{X_t}) -1 = 4$. We take $\mu = 3$. We have $\mc{L}(3)|_{{\bf P}'} \cong \mc{O}_{{\bf P}'} (4H - \sum_{i=1}^{10} D_i).$ This is a nonspecial linear system, because any 10 points on an elliptic curve impose independent conditions on quartics in ${\bf P}^2$. It follows, that the original linear system is also nonspecial.
\end{example}

\begin{example} Let $(d,n,m) = (28, 12, 8)$, expected dimension $v = 2$. We take $\mu = 9$. We have $\mc{L}(9)|_{{\bf P'}} \cong \mc{O}_{{\bf P}'} (H + \sum_{i=1}^{12} D_i)$. This is a nonspecial linear system, and so is the original one.
\end{example}

\begin{example} Let $(d,n,m) = (38,10,12)$, expected dimension $v = -1$. We take $\mu = 13$. We have $\mc{L}(13)|_{{\bf P}'} \cong \mc{O}_{{\bf P}'} (-H + \sum_{i=1}^{10} D_i)$. This is a nonspecial linear system, and so is the original one.\footnote{This example is proved in the thesis of Gimigliano \cite{Gim} by using the Horace's method (introduced in \cite{Hir1}). The original method of Ciliberto-Miranda does not handle this example (see \cite{CM2}, pp. 4048--4049).}
\end{example}

\begin{example} Let $(d,n,m) = (57,10,18)$, expected dimension $v = 0$. We take $\mu = 19$. We have $\mc{L}(19)|_{{\bf P}'} \cong \mc{O}_{{\bf P}'} (\sum_{i=1}^{10} D_i)$. This is a nonspecial linear system, and so is the original one.
\end{example}

\begin{example} Let $(d,n,m) = (174,10,55)$, expected dimension $v=-1$. In this example, our approach does not work. Indeed, to use cor. \ref{cor_cool}, we must take $\mu = 57$. But now, $\mc{L}(57)|_{{\bf P}'} \cong \mc{O}_{{\bf P}'} (3H + \sum_{i=1}^{10} 2 D_i)$, which {\it is} special! (with $h^0 = h^1 = 10$). So, the best we can say is that $h^0 (\mc{L}|_{X_t}) \leq 10$.
\end{example}

\newpage
\bibliographystyle{amsplain}

\end{document}